\documentclass[12pt]{article}

\usepackage{amsmath,amssymb,amsbsy,amsfonts,amsthm,latexsym,
                        amsopn,amstext,amsxtra,euscript,amscd,color}
\usepackage{color,xcolor}
\usepackage{latexsym}
\usepackage{cite}
\usepackage{amsfonts}
\usepackage{amsfonts,mathrsfs}
\usepackage{graphicx}
\usepackage{psfrag}
\usepackage{subfigure}
\usepackage{url}
\usepackage{stfloats}
\usepackage{amsmath}

\usepackage{hyperref}
\hypersetup{colorlinks=true}

\newtheorem{theorem}{Theorem}
\newtheorem{lemma}{Lemma}
\newtheorem{corollary}{Corollary}

\newtheorem{definition}{Definition}

\newtheorem{conjecture}{Conjecture}

\newcommand{\quash}[1]{}

\begin{document}

\title{On $q$-nearly bent Boolean functions}
\author{Zhixiong Chen$^{1}$ and Andrew Klapper$^{2}$ \\
1. Provincial Key Laboratory of Applied Mathematics,\\
 Putian University, Putian, Fujian 351100, P. R. China\\
2. Department of Computer Science,
University of Kentucky,\\ Lexington, KY
40506-0633, USA.\\
www.cs.uky.edu/$\sim$klapper.
}

\maketitle

\begin{abstract}
For each non-constant Boolean function $q$, Klapper introduced the notion of
$q$-transforms of Boolean functions. The {\em $q$-transform} of a Boolean
function $f$ is related  to the Hamming distances from $f$ to the functions
obtainable from $q$ by nonsingular linear change of basis.

In this work we discuss the existence of $q$-nearly bent functions, a new family of Boolean functions characterized by the $q$-transform.
Let $q$ be a non-affine Boolean function. We prove that any balanced Boolean functions (linear or non-linear) are $q$-nearly bent if $q$ has weight one,
 which gives a positive answer to an open question
(whether there exist non-affine $q$-nearly bent functions) proposed by Klapper. We also prove a necessary condition for checking when a function isn't $q$-nearly bent.
\end{abstract}


\section{Introduction}

Boolean functions play a central role in cryptography. For example, they appear as feedback functions in the state changes
of nonlinear feedback shift registers, as combining functions for
nonlinear combiners \cite{MS2}, and as components in S-boxes in block ciphers.

Linear or affine functions are cryptographically vulnerable and in general
should not be used.  A critical tool in the analysis of symmetric key
cryptosystems is the Walsh-Hadamard transform \cite{Ca10,CS09,X}.  It has been used, for example, to measure how closely a given Boolean function can be approximated by an affine function.  For security of  cryptographic schemes, we generally need to use functions that are `far' from affine functions, or even from simple functions (e.g., of low degree, depending on few variables, or of low rank). This reasoning has led to a new transform (see Definition \ref{q-trans}) that is related to the Hamming distance from a fixed (not necessarily affine) function \cite{Kl16}.
We (partially with other coauthors) have explored some issues concerning the
q-transform \cite{CGK,GCK,Kl16,KC}.

The topic of this work is related to the concept of {\em bent functions},  which we will review.  It is well known that for even dimension $n$,  every bent function
achieves the upper bound $2^{n-1}-2^{n/2-1}$ on the distance from the set of affine functions.
 We introduce below some basic knowledge on Boolean functions.

Let $n$ be a positive integer, let $V_n=\{0,1\}^n$, treated as row vectors, and let
${\cal B}_n=\{f:V_n\to \{0,1\}\}$, the set of Boolean functions of dimension $n$.
We denote $(0,0,\ldots,0)\in V_n$ by $0^n$. We refer the reader to Carlet's book chapter \cite{Ca10} and Cusick
and St\u{a}nic\u{a}'s monograph \cite{CS09} for background on
Boolean functions.

For $f \in {\cal B}_n$, let $wt(f) = |\{a\in V_n:f(a)=1\}|$ be the {\em
Hamming weight} of $f$.  It is the  cardinality of the support of $f$.

For $f,g \in {\cal B}_n$, let
\[W(f,g) = \sum_{a\in V_n} (-1)^{f(a)+g(a)}\in \mathbb{Z}.\]
Also let $d(f,g) = |\{a\in V_n: f(a)\ne g(a)\}| = wt(f+g)$ be the
{\em Hamming distance} from $f$ to $g$.  Then
\[
W(f,g) = 2^n-2d(f,g),
\]
which means, in particular, that $W(f,g)$ is an even number.

Let $g_v(a)=v\cdot a\in {\cal B}_n$ be linear, where $v\in V_n$.  The
\emph{Walsh-Hadamard transform coefficient} of $f\in{\cal B}_n$ at $v$
is $W(f)(v)=W(f,g_v)$.  We recall that Parseval's identity says that
\[\sum_{v\in V_n} W(f)(v)^2 = 2^{2n}.\]

A function $f \in {\cal B}_n$ with $|W(f)(v)|=2^{n/2}$ for all $v\in V_n$ is called a {\em bent function}.
For bent functions on $n$ bits to exist it is necessary that $n$ be even.
It is well known that bent functions exist for all even $n$. For example, rank $n$ quadratic functions are bent.

In recent years, Klapper generalized the Walsh-Hadamard transform to a new
transform called the $q$-transform, which is associated with a fixed non-constant
function $q\in {\cal B}_n$  \cite{Kl16}.
Let $GL_n=GL_n(\mathbb{F}_2)$ be the group of $n \times n$ invertible matrices
over $\mathbb{F}_2$. Let $N= |GL_n| = (2^n-1)(2^n-2)(2^n-4)\ldots (2^n-2^{n-1})$.  Let $q_A$
denote the function $q_A(a) = q(aA)$ for $q\in {\cal B}_n$ and $A\in GL_n$.
Let $\mathbf{0}$ denote the $n\times n$ zero matrix.

\begin{definition}\label{q-trans}
Let $q\in {\cal B}_n$.  If $f\in {\cal B}_n$ and
$A\in GL_n$, then the {\em $q$-transform coefficient of $f$ at $A$}
is $W_q(f)(A) = W(f,q_A)$.  Also let
\[W_q(f)(\mathbf{0}) = \sum_{a\in V_n} (-1)^{f(a)}\triangleq I_f,\]
where $I_f$ is referred to as the {\em imbalance} of $f$.
\end{definition}
We note that $W_q(f)(A)$ is even for all $A\in GL_n\cup \{\mathbf{0}\}$ and $I_f=0$ if and only if $f$ is balanced.

Klapper considered the statistical behavior
of the $q$-transform with respect to two probability distributions \cite{Kl16}.
For a random variable $X$ on $GL_n\cup \{\mathbf{0}\}$, let $E'[X]$
denote the expected value of $X$ with respect to the uniform
distribution on $GL_n$.  Let $E[X]$ denote the expected value of $X$ on $GL_n\cup \{\mathbf{0}\}$ with respect to $\omega$, the probability
distribution on $GL_n\cup \{\mathbf{0}\}$ defined by
\[\omega(A) = \frac{1}{N+N/(2^n-1)}= \frac{2^n-1}{2^n N}\]
for $A\in GL_n$ and
\[\omega(\mathbf{0}) = \frac{N/(2^n-1)}{N+N/(2^n-1)}= \frac{1}{2^n},\]
where $N=|GL_n|$.
Then, for balanced $q\in {\cal B}_n$, by some computations involving $\omega$ that we shall omit \cite[p. 2801]{Kl16}, we get
\begin{equation}\label{eqn-partial}
E'[W_q(f)(A)^2] = (2^{2n}-I_f^2)/(2^n-1),
\end{equation}
and so
\begin{equation}\label{eqn-Parseval}
E[W_q(f)(A)^2] = 2^n.
\end{equation}
This is a generalization of Parseval's equation.
Equation (\ref{eqn-Parseval})  leads to the notion of $q$-bent functions.

\begin{definition}(\cite{Kl16})\label{bent-def}
Let $q\in {\cal B}_n$ be balanced.  A function $f\in {\cal B}_n$ is $q$-bent if $|W_q(f)(A)|=2^{n/2}$ for all $A\in GL_n\cup \{\mathbf{0}\}$.
\end{definition}
As with bent functions, $q$-bent functions do not exist if $n$ is odd. If $q$ is
linear and non-zero, then $f$ is  $q$-bent if and only if it is bent.
However, for balanced non-affine $q$, we have proved the following.

\begin{theorem}(\cite{CGK})\label{noqbent}
Let $n\geq 4$ be even. No $q$-bent functions exist for balanced non-affine $q\in {\cal B}_n$.
\end{theorem}

Theorem \ref{noqbent} gave a positive answer to a conjecture in Klapper's
paper \cite[Conj. 3, p. 2802]{Kl16}.
So we consider generalizations of  $q$-bent functions.
Equation (\ref{eqn-partial})  leads to the notion of $q$-nearly bent functions \cite{Kl16}, where $W_q(f)(\mathbf{0}) = I_f$ is not considered.  Due to the fact that balanced
Boolean functions are especially important in cryptography, we consider the partial expected value and second
moment of the $q$-transform of such a function. From \cite[Thms.1(1) and 4]{Kl16}, the $q$-nearly bent functions are introduced in the following definition.

\begin{definition}\label{nearlybent-def}(\cite[Def.3]{Kl16})
Let $q\in {\cal B}_n$. A balanced function $f\in {\cal B}_n$ is \emph{ $q$-nearly bent} if
$$
|W_q(f)(A)|\leq \left\lceil \left( \frac{2^{2n}-I^2_q}{2^n-1}  \right)^{1/2}\right\rceil
$$
for all $A\in GL_n$, where $I_q= \sum_{a\in V_n} (-1)^{q(a)}$ as before.
\end{definition}

Note that
Definition \ref{nearlybent-def} is valid  for odd $n$ as well as even.
It follows that if $q\in {\cal B}_n$ is a quadratic function with rank $n$,
then any non-zero linear function $f$ is $q$-nearly bent. Klapper left an open question as to  whether there exist
 non-linear (balanced) $f$ so that $f$ is $q$-nearly bent for a given non-linear $q$ \cite[p.~2803]{Kl16}.

In this  work, our main contribution is to prove that there do exist non-linear (balanced) $q$-nearly  bent.
\quash{
according to the parameter $wt(q)$.
}
We organize the work as follows. In Section \ref{Necessarycondition}
we present some necessary conditions for  $q$-nearly bent functions. In Section \ref{Existence} we give a family of $q$-nearly bent functions.
In Section \ref{Nonexistence} we prove some partial results of impossibility of being $q$-nearly bent.

Throughout this work, let
$$
\rho_q=\left\lceil \left( \frac{2^{2n}-I^2_q}{2^n-1}  \right)^{1/2}\right\rceil.
$$
Since $I_{1+q}=-I_{q}$ for any $q\in {\cal B}_n$,  we see that a balanced $f$ is $(1+q)$-nearly bent if  $f$ is $q$-nearly bent.
Thus hereafter we assume that $wt(q)\leq 2^{n-1}$.

\section{Necessary conditions for $q$-nearly bent functions}\label{Necessarycondition}

In the section we show certain properties of the parameter $\rho_q$ for
any $q$-nearly bent functions. We first give two lemmas.

\begin{lemma}\label{wqf}
Let $n>2$ be any integer (odd or even) and $q\in {\cal B}_n$. For any balanced $f\in {\cal B}_n$, if $wt(q)=w\leq 2^{n-1}$, we have
for all $A\in GL_n$
$$
|W_q(f)(A)|\in \{0,4,8,\ldots,2w-4,2w\}
$$
if $w$ is even, and otherwise
$$
|W_q(f)(A)|\in \{2,6,10,\ldots,2w-4,2w\}.
$$
\end{lemma}
Proof. Since $wt(f)=2^{n-1}$ and $wt(q_A)=wt(q)=w$  for all $A\in GL_n$, we derive
$$
wt(f+q_A)\in \{2^{n-1}+w, 2^{n-1}+w-2, 2^{n-1}+w-4, \ldots, 2^{n-1}-w\}.
$$
The lemma follows from the fact that $W_q(f)(A)=2^n-2wt(f+q_A)$,
which was shown in our earlier paper \cite{KC}.        \qed

\begin{lemma}\label{average-bound}
For any $q\in {\cal B}_n$, we have
$$
(\rho_q-1)^2<\frac{2^{2n}-I^2_q}{2^n-1}\leq \rho_q^2.
$$
\end{lemma}
Proof. These bounds can be proven using the the following
$$
 \left( \frac{2^{2n}-I^2_q}{2^n-1}  \right)^{1/2}\leq \rho_q ~~~~ \text{and}  ~~~~ \left( \frac{2^{2n}-I^2_q}{2^n-1}  \right)^{1/2}+1> \rho_q.
$$
\qed

Now we study the parameter $\rho_q$.

\begin{theorem}\label{nearlybent-check-1}
Let $n>2$ be any integer (odd or even) and $q\in {\cal B}_n$. If $f\in {\cal B}_n$ is (balanced) $q$-nearly bent,
then there exists a matrix $A_0\in GL_n$ such that
$|W_q(f)(A_0)|= \rho_q$.
\end{theorem}
Proof. Suppose $|W_q(f)(A)|< \rho_q$ for all $A\in GL_n$. This implies that
$|W_q(f)(A)|^2\leq (\rho_q-1)^2$  for all $A\in GL_n$, and hence that
$E'[W_q(f)(A)^2] \leq (\rho_q-1)^2$. However, from Equation (\ref{eqn-partial})\footnote{Indeed here we use $W_q(f)(A)=W_f(q)(A^{-1})$ by Thm. 1(1) in \cite{Kl16}.}
 and Lemma \ref{average-bound}, we see that $(\rho_q-1)^2<E'[W_q(f)(A)^2]\leq \rho_q^2$,
a contradiction. \qed

\begin{theorem}\label{nearlybent-check-2}
Let $n>2$ be any integer (odd or even) and $q\in {\cal B}_n$. If there exists a (balanced) $q$-nearly bent function in ${\cal B}_n$,
then,  either  $\rho_q \equiv 0 \pmod 4$ when $wt(q)$ is even
or  $\rho_q \equiv 2 \pmod 4$ when $wt(q)$ is odd.
\end{theorem}
Proof. Suppose that $f$ is $q$-nearly bent. Then there exists a matrix $A_0\in GL_n$ such that
$|W_q(f)(A_0)|= \rho_q$ by Theorem \ref{nearlybent-check-1}. The statements can be proved by Lemma \ref{wqf}.
 \qed

Theorem \ref{nearlybent-check-2} can be used to check whether there exist
$q$-nearly bent functions for certain $q$.
We prove the non-existence of $q$-nearly bent functions for certain $q$ in Section \ref{Nonexistence}.
But in the coming section, we present a family of $q$-nearly bent functions.

\section{Existence of $q$-nearly bent functions}\label{Existence}

Since linear functions are $q$-nearly bent for any bent function $q$, Klapper asked
whether there exist non-linear $q$ and non-linear balanced $f$ so that $f$ is  $q$-nearly bent  \cite[p.~2803]{Kl16}.
Theorem \ref{nearly-wt1} below answers this question.

\begin{theorem}\label{nearly-wt1}
Let $n>2$ be any integer (odd or even) and $q\in {\cal B}_n$. If $wt(q)=1$, then any balanced $f$ (linear or non-linear) is $q$-nearly bent.
\end{theorem}
Proof. By Lemma \ref{wqf} we see that $W_q(f)(A)\in \{-2,2\}$ for all $A\in GL_n$. Now using the fact that $I_q=2^n-2wt(q)=2^n-2$, we can prove
$$
|W_q(f)(A)|=2=\left( \frac{2^{2n}-(2^n-2)^2}{2^n-1}  \right)^{1/2} = \rho_q.
$$
So by Definition \ref{nearlybent-def}, $f$ is $q$-nearly bent. \qed

It is natural to consider the cases of other values of $wt(q)$. In the following section, we find different phenomena.

\section{Impossibility  of being $q$-nearly bent}\label{Nonexistence}

In this section we derive conditions under which no $q$-nearly bent
functions can exist.  These conditions are expressed in terms of $wt(q)$
and $\rho_q$.  They apply when the weight of $q$ is low.

\begin{theorem}\label{nearly-wt2,3}
Let $n>2$ be an integer (odd or even) and $q\in {\cal B}_n$. If $wt(q)\in \{2,3\}$, then there are no non-linear $q$-nearly bent functions.
\end{theorem}
Proof. Using $I_q=2^n-2wt(q)$, we compute
\[
 \rho_q=\left\{
\begin{array}{ll}
3, & \mathrm{if}\,\ wt(q)=2 \,\  \mathrm{and}\,\ n\geq 3,  \\
3, & \mathrm{if}\,\ wt(q)=3  \,\  \mathrm{and}\,\ n=3,  \\
4, & \mathrm{if}\,\ wt(q)=3  \,\  \mathrm{and}\,\ n>3.
\end{array}
\right.
\]
Then the proof is finished using Theorem \ref{nearlybent-check-2}. \qed

Table \ref{table-experiment} below lists some experimental results.
Some of these results can be checked using Theorem
\ref{nearlybent-check-2}, but others cannot.

\begin{table}[hb]
\centering
\begin{tabular}{cccl}
\hline\noalign{\smallskip}
$n$ & $wt(q)$ & $\rho_q$        &  Do $q$-nearly bent functions exist?  \\
 \noalign{\smallskip} \hline \noalign{\smallskip}
$n>2$ & $4$  &  $\rho_q=4$           &  No (Proved in Thm.\ref{nearlybent-wt4})                  \\
$n>3$ & $5$  & $\rho_q\in\{4,5\}$    & No (by Thm.\ref{nearlybent-check-2})                \\
$n=4$ & $6$  & $\rho_q=4$           & No (Proved in Thm.\ref{nearlybent-largewtq})                 \\
$n>4$ & $6$  & $\rho_q=5$           &  No (by Thm.\ref{nearlybent-check-2})                \\
$n=4$ & $7$   & $\rho_q=2$          &  No (Proved in Thm.\ref{nearlybent-largewtq})                 \\
$n=5$ & $7$  & $\rho_q=5$           & No  (by Thm.\ref{nearlybent-check-2})               \\
$n>5$ & $7$  & $\rho_q=6$           & ?                 \\
\hline
\end{tabular}
\caption{Examples of $q$ with small $wt(q)$}\label{table-experiment}
\end{table}

We see that for $wt(q)=4$, Theorem \ref{nearlybent-check-2} is of no use to us.
So we need  further analysis.

\begin{theorem}\label{nearlybent-wt4}
Let $n>2$ be any integer (odd or even) and $q\in {\cal B}_n$. If $wt(q)=4$, then no non-linear $q$-nearly bent function exists.
\end{theorem}
Proof. We first prove the case when $n=3$. Let $q\in {\cal B}_3$ have $wt(q)=4$.
In this case  $q$ is balanced. We only need to deal with the balanced functions with value $0$ at $0^3=(0,0,0)$.
There are only $\binom{7}{4}=35$ balanced functions (with value $0$ at $0^3$), which include 7 of degree 1 and 28 of degree 2. Note that there are no balanced functions of degree 3
\footnote{There are no balanced functions with degree $n$ in ${\cal B}_n$.}.

On the other hand, the number of $3 \times 3$ invertible matrices over $\mathbb{F}_2$ is $|GL_3|=(8-1)(8-2)(8-4)=168$.
Consider the balanced function $q(x)=q(x_1,x_2,x_3)=x_1x_2+x_3\in {\cal B}_3$. There are exactly six $A\in {\cal B}_3$'s such that $q(xA)=q(x)$.  They are
$$
\left[
\begin{array}{ccc}
1 & 0 & 0\\
0 & 1 & 0\\
0 & 0 & 1
\end{array}
\right],
\left[
\begin{array}{ccc}
0 & 1 & 0\\
1 & 0 & 0\\
0 & 0 & 1
\end{array}
\right],
\left[
\begin{array}{ccc}
1 & 1 & 1\\
0 & 1 & 0\\
0 & 0 & 1
\end{array}
\right],
$$

$$
\left[
\begin{array}{ccc}
1 & 1 & 1\\
1 & 0 & 0\\
0 & 0 & 1
\end{array}
\right],
\left[
\begin{array}{ccc}
0 & 1 & 0\\
1 & 1 & 1\\
0 & 0 & 1
\end{array}
\right], \text{ and }
\left[
\begin{array}{ccc}
1 & 0 & 0\\
1 & 1 & 1\\
0 & 0 & 1
\end{array}
\right].
$$
These six matrices form $q$'s stabilizer subgroup
under the action of $GL_3$ on ${\cal B}_3$.
Hence for all $B\in GL_3$, there are  $168\div 6=28$ distinct $q_B(x)$'s.
This implies that all balanced functions of degree 2 are of the form
$q_B(x)$ for some $B\in GL_3$.

Thus for any balanced $f$ of degree 2, we have $f(x)=q(xA_0)$ for some $A_0\in GL_3$. Hence
$|W_q(f)(A_0)|=8>\rho_q=4$. That is, no non-linear $q$-nearly bent function exists for balanced $q\in {\cal B}_3$ of degree 2.

Now we turn to the case when $n\geq 4$. Since $wt(q)=4$ we may assume
that $q(0^n)=0$ and select $\alpha_i\in V_n$ for $1\leq i\leq 4$ such that
$q(\alpha_i)=1$, and $q(\beta)=0$ for $\beta\not\in \{\alpha_i:i=1,\cdots,4\}$.
We will show that for any non-linear balanced $f$ there exists an $M\in GL_n$ such that $|W_q(f)(M)|= 8$. This implies that $|W_q(f)(M)|>\rho_q=4$ and hence $f$  isn't $q$-nearly bent.
If $f(\alpha_i)=0$ for all $i\in \{1,2,3,4\}$ or $f(\alpha_i)=1$ for
all $i\in \{1,2,3,4\}$, we get immediately that $|W_q(f)(E)|= 8$, where $E$ is the identity matrix.  Otherwise, we give a proof in two cases. We assume
below that $f$ is non-linear and balanced, and that $f(0^n)=0$.

Suppose $\alpha_1,\alpha_2,\alpha_3,\alpha_4$ are linearly independent.
The support of $f$ has cardinality $2^{n-1}$ and is not contained in a
proper subspace of $V_n$ (otherwise the
support would equal the subspace, which would imply that $f$ is affine).
We select four row vectors $\beta_1,\beta_2,\beta_3,\beta_4 \in \text{supp}(f)$
that are linearly independent.  Hence there exists an $M\in GL_n$  such that
$\alpha_i=\beta_i M$ for $1\leq i\leq 4$. Then for $1\leq i\leq 4$ we have
$$
q(\beta_i M)=q(\alpha_i)=1 ~~~\text{and}~~~ f(\beta_i)=1.
$$
So we get
$$
|W_q(f)(M)|=\left|\sum_{a\in V_n} (-1)^{f(a)+q(aM)}\right| =8.
$$

Now suppose $\alpha_1,\alpha_2,\alpha_3,\alpha_4$ are linearly dependent.
Without loss of generality, we suppose  that $\alpha_1,\alpha_2,\alpha_3$ are linearly independent and $\alpha_4$ is of the form $\alpha_4=\alpha_1+\alpha_2+\alpha_3$ or $\alpha_4=\alpha_1+\alpha_2$.

First, if $\alpha_4=\alpha_1+\alpha_2+\alpha_3$, then there are $2^{n-1}(2^{n-1}-1)$ vectors of the form $\gamma_i+\gamma_j$ with
$\gamma_i\ne \gamma_j\in \text{supp}(f)$.  We can select four pairwise distinct
vectors $\gamma_1,\gamma_2,\gamma_3,\gamma_4\in \text{supp}(f)$ such
that  $\gamma_1+\gamma_2=\gamma_3+\gamma_4$.  Hence there exists an $M\in GL_n$  such that
$\alpha_i=\gamma_i M$ for $1\leq i\leq 3$, from which we get
$\alpha_4=\gamma_4 M$ and
$$
q(\gamma_i M)=q(\alpha_i)=1 ~~~\text{and}~~~ f(\gamma_i)=1
$$
for $1\leq i\leq 4$. So we get $|W_q(f)(M)|=8$.

Second, we suppose $\alpha_4=\alpha_1+\alpha_2$. If there exist
pairwise distinct $\gamma_1,\gamma_2,\gamma_3\in \text{supp}(f)$ such that  $\gamma_3\neq \gamma_1+\gamma_2$,
then $\gamma_1,\gamma_2,\gamma_3$ are linearly independent.
We write $\gamma_4\triangleq \gamma_1+\gamma_2$. If $\gamma_4\in \text{supp}(f)$, and then
we can construct $M\in GL_n$  such that
$\alpha_i=\gamma_i M$ for $1\leq i\leq 3$ and calculate $|W_q(f)(M)|=8$.
If $\gamma_4\not\in \text{supp}(f)$, that is, if any vector $\gamma_i+\gamma_j\in V_n\setminus\text{supp}(f)$ for distinct $\gamma_i,\gamma_j\in \text{supp}(f)$,
then we choose pairwise distinct $\gamma_1,\gamma_2,\gamma_3\in \text{supp}(f)$ and set
$$
\eta_1=\gamma_1+\gamma_2,~~\eta_2=\gamma_1+\gamma_3,~~\eta_4=\gamma_2+\gamma_3.
$$
Then $\eta_1,\eta_2,\eta_4 \in V_n\setminus\text{supp}(f)$. We further select $0^n\neq \eta_3\in V_n\setminus\text{supp}(f)$ such that $\eta_3\neq \eta_4 =\eta_1+\eta_2$. Hence we can construct $M\in GL_n$  such that
$\alpha_i=\eta_i M$ for $1\leq i\leq 3$, from which we get $\alpha_4=\eta_4 M$ and
$$
q(\eta_i M)=q(\alpha_i)=1 ~~~\text{but}~~~ f(\eta_i)=0
$$
for $1\leq i\leq 4$. So we get $|W_q(f)(M)|=8$.  This completes the proof. \qed\\

We hope that  the method used  in Theorem \ref{nearlybent-wt4} can be useful for checking  $q$ with other $wt(q)$.
For example, if $4<wt(q)\leq n$ and all vectors in $\mathrm{supp}(q)$ are  linearly independent, then
any non-linear balanced function isn't $q$-nearly bent. However, if the vectors in $\mathrm{supp}(q)$ are  linearly dependent,
the proof seems to be more complicated.

We remark that the $wt(q)$ considered above is small, so we prove a result for large $wt(q)$.

\begin{theorem}\label{nearlybent-largewtq}
Let $n>2$ be  even and $q\in {\cal B}_n$.
If
\begin{equation}\label{Eqn:nearlybent-largewtq}
2^{n-1}-2^{n/2-1}\le wt(q)\leq 2^{n-1},
\end{equation}
then no non-linear $q$-nearly bent function exists.
\end{theorem}

Proof. Suppose that non-linear, balanced $f$ is $q$-nearly bent.
We have $|\{x: q(x)=1\}| = wt(q)$, and $|\{x: q(x)=0\}| = 2^n -wt(q)$.
Thus $I_q = |\{x:q(x)=0\}|- |\{x: q(x)=1\}| = 2^n-2wt(q)$.
\quash{From equation
(\ref{Eqn:nearlybent-largewtq}) and the fact that $0\le I^2_q\leq 2^n$, we get
the equality $I_q=2^n-2wt(q)$.
}
We have
$$
\begin{array}{ll}\displaystyle
2^{n}+1< \frac{2^{2n}-I^2_q}{2^n-1}  <2^{n}+2,  & \text{if}~~~ I^2_q=0,\\[12pt]
\displaystyle 2^{n}< \frac{2^{2n}-I^2_q}{2^n-1}  <2^{n}+1,    &
              \text{if}~~~ 0< I^2_q< 2^n, \text{and}\\[12pt]
\displaystyle \frac{2^{2n}-I^2_q}{2^n-1} =2^{n},            & \text{if}~~~ I^2_q= 2^n.
\end{array}
$$
Note that $\sqrt{2^{n}+2}<2^{n/2}+1$. Thus for $0\le I^2_q< 2^n$ we have
$$
2^{n/2}<\rho_q\leq \lceil \sqrt{2^{n}+2}\rceil \leq 2^{n/2}+1,
$$
and hence $\rho_q= 2^{n/2}+1$ (an odd number), which is impossible by Theorem \ref{nearlybent-check-2}.

Now for $I^2_q= 2^n$ and $A\in GL_n$ we have by Definition \ref{nearlybent-def}
$$
|W_q(f)(A)|\leq 2^{n/2}.
$$
Then we see that $|W_f(q)(A)|=|W_q(f)(A^{-1})|\leq 2^{n/2}$ for all $A\in GL_n$ and $|W_f(q)(\mathbf{0})|=I_q= 2^{n/2}$.
Equation (\ref{eqn-Parseval}) implies that $|W_f(q)(A)|= 2^{n/2}$ for all $A\in GL_n\cup\{\mathbf{0}\}$ and hence $q$ is $f$-bent,
which contradicts Theorem \ref{noqbent}.  This completes the proof. \qed\\

We have immediately the following corollary, which says Proposition 5 in \cite{Kl16} is vacuous for even $n$.

\begin{corollary}
Let $n>2$ be even and $q\in {\cal B}_n$. If $q$ is balanced (so $wt(q)=2^{n-1}$), then no non-linear (balanced) $q$-nearly bent function exists.
\end{corollary}

Unfortunately, we cannot prove a similar result based on Theorem \ref{nearlybent-largewtq} for odd $n\geq 3$. For $n=3$, the first part of the proof of Theorem \ref{nearlybent-wt4}
says that any non-linear balanced $f\in {\cal B}_3$  of degree 2  is not $q$-nearly bent for non-linear balanced $q\in {\cal B}_3$. For odd $n>3$, we haven't any results,
although sometimes it can be checked by Theorem \ref{nearlybent-check-2}.
We list blow some experimental data. We cannot prove the cases when $n=7$ or $n=9$.

\begin{center}
\begin{tabular}{ccclc}
\hline\noalign{\smallskip}
&$n$ &  $\rho_q$        &  Do $q$-nearly bent functions exist?  &\\
 \noalign{\smallskip} \hline \noalign{\smallskip}
&3 &  $4$         & No (see Thm.\ref{nearlybent-wt4})     & \\
&5 &  $6$         & No (by Thm.\ref{nearlybent-check-2})     &\\
&7 &  $12$        & ?    &\\
&9 &  $24$        & ?     & \\
&11 &  $46$        & No (by Thm.\ref{nearlybent-check-2})         &  \\
&13 &  $91$        & No  (by Thm.\ref{nearlybent-check-2})       &     \\
&15 &  $182$       & No  (by Thm.\ref{nearlybent-check-2})         &     \\
&17 &  $363$       & No  (by Thm.\ref{nearlybent-check-2})      &         \\
\hline
\end{tabular}
\\ Table 2. Examples of $q$ with $wt(q)=2^{n-1}$
\end{center}

For $q(x)=q(x_1,x_2,\ldots,x_n)=x_1x_2\in {\cal B}_n$, we see that $wt(q)=2^{n-2}(<2^{n-1}-2^{n/2-1})$ and
check that there are no $q$-nearly bent functions when $n=4$ by Theorem \ref{nearlybent-wt4} and when $5\leq n\leq 8$ by Theorem \ref{nearlybent-check-2} (by checking $\rho_q$).
We do not know whether $q$-nearly bent functions exist if $n=9$, in which case $\rho_q=20$.

The results and examples seem to indicate that there are no non-linear $q$-nearly bent functions when $1<wt(q)\leq 2^{n-1}$, so we make this
a conjecture.

\begin{conjecture}\label{conj}
Let $n>2$ be any integer (odd or even) and $q\in {\cal B}_n$. If $1< wt(q)\leq 2^{n-1}$,
 then no non-linear (balanced) $q$-nearly bent function exists.
\end{conjecture}

\section{Conclusions}\label{Conclusion}

In this work, we have studied the existence and non-existence of $q$-nearly bent functions with respect to $q$-transform for  a non-affine Boolean function $q$.
We have shown that all balanced functions are $q$-nearly bent
if $q$ has weight one, which confirms the existence of $q$-nearly bent functions.
We have also pointed out a necessary condition for $q$-nearly bent
functions to exist. That is, 
 if a $q$-nearly bent function exists then either  $\rho_q \not\equiv 0 \pmod 4$ when $wt(q)$ is even
or  $\rho_q \not\equiv 2 \pmod 4$ when $wt(q)$ is odd. The statement can be used to determine the non-existence of $q$-nearly bent functions.
Results and examples lead us to the conjecture that there are no non-linear $q$-nearly bent functions when $1<wt(q)\leq 2^{n-1}$.

To weaken the notion of q-bentness, it is natural
 to generalize plateaued functions \cite{HLL,ZZ},
  we would say $f\in {\cal{B}}_n$ is \emph{$q$-plateaued} if all its $q$-transform coefficients are in $\{0, \pm\lambda\}$
for some positive integer $\lambda$.
When $n=3$ and $q(x_1,x_2,x_3)=x_1x_2+x_3$, we can check that any $f\in {\cal{B}}_3$ with $wt(f)=2$ is $q$-plateaued. In this case,
$\lambda=4$. Then we ask,  for $n>3$ and a given non-affine $q$, do there exist non-affine $q$-plateaued functions?

\section*{Acknowledgments}

Z. Chen was partially supported by the National Natural Science
Foundation of China under grant No.~61772292, by the Provincial Natural Science
Foundation of Fujian under grant No.~2018J01425 and by the Program for Innovative Research Team in Science and Technology in Fujian Province University under grant No.~2018-49.

A. Klapper was partially supported by the National Science Foundation
under Grant No.~CNS-1420227. Any opinions, findings, and conclusions or recommendations expressed in this material are those of the authors and do not necessarily reflect the views of the National Science Foundation.

\end{document}